%% file: 9-13-040-Quasi.tex
\documentclass{article}

\usepackage{graphicx}
\usepackage{amsfonts}
\usepackage{amsmath}
\usepackage{amsthm}
\usepackage{epsfig,color}

\newtheorem{Thm}{Theorem}

\newtheorem{Lem}{Lemma}

\newtheorem{Def}{Definition}

\newcommand{\q}{\mathbb{Q}}

\begin{document}
\title{Quasiinvariants of $S_3$}
\author{Jason Bandlow, Gregg Musiker}
\date{September 13, 2004}
\maketitle

\begin{abstract} \noindent Let $s_{ij}$ represent a transposition in
$S_n$.  A polynomial $P$ in $\q[X_n]$ is said to be
$m$-quasiinvariant with respect to $S_n$ if $(x_i-x_j)^{2m+1}$
divides $(1-s_{ij})P$ for all $1 \leq i, j \leq n.$ We call the
ring of $m$-quasiinvariants, $QI_m[X_n]$.  We describe a method
for constructing a basis for the quotient
$QI_m[X_3]/(e_1,e_2,e_3)$.  This leads to the evaluation of
certain binomial determinants that are interesting in their own
right.

\end{abstract}

\noindent The symmetric group $S_n$ acts on the ring of
polynomials $\q[X_n]$ by permuting indices.  That is for any
permutation $\sigma \in S_n$
\begin{align*}
\sigma  P(x_1,\dots, x_n) = P(x_{\sigma(1)},\dots, x_{\sigma(n)}).
\end{align*}

\noindent A polynomial $P$ is said to be $S_n$-invariant or
symmetric if and only if $\sigma(P) = P \mathrm{~~for~all~~}
\sigma \in S_n.$  The fundamental theorem of symmetric functions
\cite[p. 292]{EC2} states that any invariant of $S_n$ can be
written as a polynomial in $\{e_1,e_2,\dots, e_n\}$ where
\begin{align*}
e_k = \sum_{1 \leq i_1 < \dots < i_k \leq n} x_{i_1}x_{i_2}\cdots
x_{i_k}. \end{align*} For $S_3$ we have
\begin{align*}
e_1 &= x_1 + x_2+x_3 \\
e_2 &= x_1x_2 + x_1x_3 + x_2x_3 \\
e_3 &= x_1x_2x_3.
\end{align*}

\noindent A generalization of invariance known as
``\emph{quasiinvariance}" has been studied in the recent
literature \cite{FV,FV2,EG}. In the rest of this paper we will use
the notation $s_{ij}$ to denote the transposition $(i,j)$ and will
let $QI_m$ denote $QI_m[X_n]$ for convenience.

\begin{Def}
A polynomial $P$ is $m$-quasiinvariant if and only if $(1-s_{ij})P$ is divisible by
$(x_i-x_j)^{2m+1}$ for all pairs $1\leq i < j \leq n$.
\end{Def}

\noindent This definition is not vacuous because $(1-s_{ij})P$ is
antisymmetric with respect to the transposition $s_{ij}$ thus
setting $x_i = x_j$ will yield zero.  Hence $(x_i-x_j)$ divides
$(1-s_{ij})P$ and the antisymmetry forces an odd power of
$(x_i-x_j)$ to divide it.  We should note that an analogous
condition defines $m$-quasiinvariance for any Coxeter group.  In
the general definition, the linear forms giving the equations of
the reflecting hyperplanes play the role of the differences
$x_i-x_j$.

\vspace{1.5em} \noindent It is easily seen that the divided
difference operator $\Delta_{ij} = \frac{1 -s_{ij}}{ x_i-x_j }$ is
a twisted derivation \cite[pp. 192-194]{Kane} which means that

$$\Delta_{ij}(Q_1Q_2) = \Delta_{ij}(Q_1)Q_2 + s_{ij}(Q_1)\Delta_{ij}(Q_2).$$

\noindent Thus if $\Delta_{ij}(Q_1)$ and $\Delta_{ij}(Q_2)$ are
both divisible by $(x_i-x_j)^{2m}$ then so is
$\Delta_{ij}(Q_1Q_2)$.  The operator $(1-s_{ij})$ is also linear
which means that each $QI_m$ is a ring.  Furthermore $(1-s_{ij})P$
will be divisible by $(x_i-x_j)^{2m+1}$ for arbitrarily large $m$
if and only if $(1-s_{ij})P = 0$ which means that all $QI_m$
contain $\Lambda_n$ (the ring of symmetric polynomials).  We thus
have the inclusions

$$\q[x_1,\dots, x_n] = QI_0 \supset QI_1 \supset QI_2 \supset \cdots \supset QI_\infty = \Lambda_n$$

\vspace{1.5em} \noindent A classic result states that $\q[x_1,\dots,x_n]$ is a free module of rank $n!$ over the ideal
$(e_1,\dots, e_n)$.  Furthermore, the action on the quotient precisely gives the regular representation of $S_n$
\cite[p. 247]{Kane}.

\vspace{1.5em} \noindent This means that there exists a basis of $n!$ polynomials
$\{\eta_1,\dots,\eta_{n!}\}$ such that any $n$-variable polynomial can be written as a unique linear
combination
$$\sum_{i=1}^{n!} A_i\eta_i$$ where the $A_i$'s are symmetric polynomials.  For example, any
polynomial in $\q[x_1,x_2,x_3]$ can be written uniquely as
$$A_1 + A_2 x_2 + A_3 x_3 + A_4 x_2x_3 + A_5 x_3^2 + A_6 x_2x_3^2$$

\noindent where $A_1,\dots, A_6$ are symmetric polynomials.  The
polynomial ring can be thought of as the ring of
$0$-quasiinvariants and recently \cite{EG}, an analogous result
has been proven for the rings of $m$-quasiinvariants for $m > 0$.
Namely, any element of $QI_m$ can be written uniquely as a sum

$$\sum_{i=1}^{n!} A_i(e_1,\dots, e_n) \cdot \eta_i$$

\noindent where the $A_i$'s are polynomials and the $\eta_i$'s are
elements of $QI_m$.

\vspace{1.5em} \noindent These $\eta_i$'s are therefore a basis
for $QI_m \big/ \big\langle (e_1,e_2, \dots e_n) \big\rangle$, a
space which has been shown \cite{FV2} to have the following
Hilbert series:

$$\sum_{i=1}^{n!} q^{degree(\eta_i)} =  \sum_{T \in ST(n)} q^{m \left( \binom{n}{2} - content(\lambda(T) )\right)+cocharge(T)}$$

\vspace{1.5em} \noindent In the case that $n=3$, this gives that
the Hilbert series of $QI_m \big/ \big\langle (e_1,e_2,e_3)
\big\rangle$ ($QI_m$ will always signify $QI_m[X_3]$ from here on
out) is

\begin{align} \label{hilb} q^0 + 2q^{3m+1} + 2q^{3m+2} + q^{6m+3}. \end{align}

\noindent Note also by the respective degrees of $e_1, e_2,$ and
$e_3$ that the Hilbert series of $QI_m$ is

\begin{align} \label{hilb2} {q^0 + 2q^{3m+1} + 2q^{3m+2} + q^{6m+3} \over (1-q)(1-q^2)(1-q^3)}. \end{align}

\noindent It is easily shown that the Vandermonde determinant
$\Delta(x) = (x_1 - x_2)(x_1-x_3)(x_2-x_3)$ raised to the power
$2m+1$ accounts for the term $q^{6m+3}$ and clearly the constants
account for $q^0$. So the interesting problem arises to construct
the four $m$-quasiinvariants that account for the terms $2q^{3m+1}
+ 2q^{3m+2}$.  The explicit construction of these four
$m$-quasiinvariants is the goal and motivating force which led to
the results of this paper.  It developed that this construction
required the evaluation of two binomial determinants which are
interesting in their own right and deserve a special mention here.
The two resulting identities may be stated as follows.

\begin{Thm} \label{specthm}
\begin{align} \label{specdet}
\det \Bigg| \binom{C+\alpha i}{E+\beta j} - \binom{D-\alpha
i}{E+\beta j} \Bigg|_{i,j=1}^k =
\frac{\binom{C+D}{E+\beta}\binom{C+D}{E+2\beta} \dots
\binom{C+D}{E+n\beta}} {\binom{C+D}{C+
\alpha}\binom{C+D}{C+2\alpha} \dots \binom{C+D}{C+n\alpha}} \cdot
| \mathcal{F} |
\end{align}
where $\mathcal{F}$ denotes the collection of $k$-tuples of
non-intersecting lattice paths respectively joining the points
$$\{(D-k\alpha,D-k\alpha),(D-(k-1)\alpha,D-(k-1)\alpha),\dots,(D-\alpha,D-\alpha)\}$$
to
$$\{(0,C+D-E -k\beta),(0,C+D-E -(k-1)\beta),\dots,(0,C+D-E-\beta)\}$$ and
throughout remaining strictly below the line $y=-x+C+D$.
\end{Thm}

\noindent It is also worthy of notice the fact that the entries of
the determinant in (\ref{specdet}) are differences of binomial
coefficients where the tops are different and the bottoms are the
same.  A literature search found no determinant results covering
this particular case.  Nevertheless, a manipulation suggested by
an argument of Gessel and Viennot in \cite{Gessel} enabled us to
derive Theorem \ref{specthm} from the following general result:

\begin{Thm} \label{Gendet} For any integers $a,b,c,d,e$, the determinant
\[
\det \Bigg| \binom{a+bi}{c+dj} - \binom{a+bi}{e-dj}
\Bigg|_{i,j=1}^n
\]
is the number of families of non-intersecting lattice paths with
NORTH and WEST steps, respectively joining the points
$$\{(c+d,c+d), (c+2d,c+2d), \dots (c+nd,c+nd) \}$$ to
$$\{(0,a+b),(0,a+2b),\dots,(0,a+nb)\}$$ and throughout avoiding the line $y=-x + (c+e)$.
\end{Thm}

\noindent Our main result is that a basis for the quotient of the
$m$-quasiinvariants of $S_3$ can be found by computing the
1-dimensional null space of particular matrices. The non-vanishing
of the determinant (\ref{specdet}) provides the crucial step in
proving the null space in question is indeed 1-dimensional.

\vspace{1.5em}\noindent Our presentation is divided into four
parts.  In the first part we show (non-constructively) that
quasiinvariants of a certain nice form exist.  In the second part,
we find a system of equations that the coefficients of these
quasiinvariants must satisfy.  In the third part, we show that we
can solve this system by computing a 1-dimensional null space.  In
the final part we complete the construction, and prove the
elements we've constructed complete a basis for the quotient.

\vspace{1.5em}\noindent We should mention that Feigin and Veselov
in \cite{FV} have given explicit module bases for the
$m$-quasiinvariants of all Dihedral groups $D_n$.  But so far
there are no other Coxeter groups for which explicit constructions
have been given.  The Feigin-Veselov construction is based on
complex number techniques that are very suitable in the dihedral
case.  Although $D_3$ $m$-quasiinvariants can be easily converted
into $S_3$ $m$-quasiinvariants, our work efforts have been guided
by the need of developing methods that can be extended to the
general case.  Our results may be taken as an instance of such
methods.  Extensions of the present construction to $S_n$ will be
the topic of a forthcoming publication.

\section{Quasiinvariants with a nice form}

\noindent We begin by defining the following elements of the group
algebra of $S_3$:
\begin{align*}
[S_3] &= \frac{1}{6}\sum_{\sigma \in S_3} \sigma,&
[S_3]' &= \frac{1}{6}\sum_{\sigma \in S_3} \mathrm{sgn}(\sigma) \sigma\\
\pi_1 &= \frac{1}{3}(1+s_{23})(1-s_{12}), & \pi_2
&=\frac{1}{3}(1+s_{12})(1-s_{23})
\end{align*}

\noindent These defined, the following identities are easily
verified:
\begin{eqnarray}
(\pi_1)^2=\pi_1, (\pi_2)^2=\pi_2 \label{idemp1}\\
\lbrack S_3]' \pi_1 = \pi_1 \pi_2 = \pi_2 \pi_1= 0 \label{idemp2}\\
\lbrack S_3] + \pi_1 + \pi_2 + [S_3]' = 1 \label{ident} \\
s_{23}\pi_1 = \pi_1 \label{sym12}\\
\pi_2 s_{12} \pi_1 = - s_{13} \pi_1 \label{swap}
\end{eqnarray}

\noindent We now show that there exist quasiinvariants satisfying
certain symmetry and independence conditions.

\begin{Lem} \label{linind} For all $m \ge 0$, there exist
non-symmetric $m$-quasiinvariants $A_1, A_2$ of degrees $3m+1,
3m+2$, respectively, such that $s_{23} (A_i) = A_i$ and in the
quotient $QI_m \big/ \big\langle (e_1,e_2, e_3) \big\rangle$, the
image of $A_i$ and $s_{12} (A_i)$ are linearly independent.
Further all four of these will be independent of $\Delta^{2m+1}(x)$.
\end{Lem}

\begin{proof}
\noindent It is easy to see that the image of $[S_3]$ in the
quotient is the constant terms.  We also note that any polynomial
in the image of $[S_3]'$ is alternating and any alternating
$m$-quasiinvariant must be divisible by $\Delta^{2m}(x)$, which
has degree $6m$. Thus, from the Hilbert series (\ref{hilb}), there
must exist quasiinvariants $B_i$ of degree $3m+i$, ($i \in
\{1,2\}$) such that if we apply equation~(\ref{ident}) to $B_i$ we
have
\begin{align}
\pi_1(B_i) + \pi_2(B_i) \ne 0.
\end{align}
Assume without loss that $\pi_1(B_i) \ne 0$, and set
\begin{align} A_i = \pi_1(B_i).
\end{align}
Equation (\ref{sym12}) immediately gives that $s_{23}(A_i) = A_i$.
Now suppose we had symmetric functions $S,T$ such that
\begin{align} \label{depend} SA_i + T (s_{12}A_i) = 0
\end{align}
Applying $\pi_2$ to this gives (by (\ref{idemp2}) and
(\ref{swap})):
\begin{align}
T \pi_2 s_{12} \pi_1 B_i = 0\\
-T s_{13} \pi_1 B_i = 0\\
-T s_{13} A_i = 0
\end{align}
Since $A_i$ was assumed to be non-zero, this gives $T=0$ and
(\ref{depend}) gives $S=0$. Now assume there was a nontrivial
relationship between these and $\Delta^{2m+1}(x)$
\begin{align}
c_1A_1 +  c_2(s_{12}A_1) + c_3A_2 +  c_4(s_{12}A_2) + c_5\Delta^{2m+1}(x)
= 0
\end{align}
Applying $[S_3]'$ gives (by (\ref{idemp2}))
\begin{align}\label{last}
c_2[S_3]'(s_{12}A_1) + c_4[S_3]'(s_{12}A_2) + c_5\Delta^{2m+1}
(x) = 0
\end{align}
But $[S_3]'s_{12} = s_{12}[S_3]'$ and $[S_3]'A_i = 0$ so
(\ref{last}) gives $c_5 =0$.
\end{proof}

\noindent Since $s_{23}(A_i)=A_i$, $A_i$ is symmetric with respect
to $x_2$ and $x_3$. This means that we can write the
$m$-quasiinvariants $A_1$ and $A_2$ as

\begin{align} \label{thesum} A_1 = \sum_{\substack{0 \leq i \leq j \leq i+j \leq d }} C_{[i,j]} x_1^{d-i-j} m_{[i,j]}(x_2,x_3).
\end{align}
and
\begin{align} A_2 = \sum_{\substack{0 \leq i \leq j \leq i+j \leq d }} \tilde{C}_{[i,j]} x_1^{d-i-j} m_{[i,j]}(x_2,x_3).
\end{align}

\noindent for $d = 3m+1$ or $3m+2$, respectively.  In fact we can
make the following stronger statement about the form of the $A_i$:

\begin{Lem} \label{betterform} There exist $m$-quasiinvariants $A_1$ and
$A_2$, satisfying the conditions of Lemma~\ref{linind}, of the
form
\begin{align*}
A_1 = \sum_{0 \leq i \leq j \leq m} C_{[i,j]} x_1^{3m+1-i-j}
m_{[i,j]}(x_2,x_3) \
\end{align*}
and \begin{align*} A_2 = \sum_{0 \leq i \leq j \leq m+1}
\tilde{C}_{[i,j]} x_1^{3m+2-i-j} m_{[i,j]}(x_2,x_3).
\end{align*}
\end{Lem}

\begin{proof} We first prove this result for $A_1$.  By grouping together monomials with similar exponent sequences,
we can rewrite the above sum (\ref{thesum}) as
\begin{align*}
&\sum_{\stackrel{0 \leq i < j < k}{i+j+k ~=~ 3m+1}}  \bigg(
C_{[i,j]} (x_1^kx_2^ix_3^j  + x_1^kx_2^jx_3^i) ~+~ C_{[j,k]}
(x_1^ix_2^jx_3^k  + x_1^ix_2^kx_3^j) \\
& ~+~ C_{[i,k]} (x_1^jx_2^kx_3^i  + x_1^jx_2^ix_3^k) \bigg)  +
\sum_{\stackrel{0 \leq i ,  j}{2i+j ~=~ 3m+1}}  C_{[i,j]}
(x_1^ix_2^ix_3^j ~+~ x_1^ix_2^jx_3^i) \\
& ~+~ C_{[i,i]}x_1^jx_2^ix_3^i.
\end{align*}

\noindent Using this decomposition, we find that $(1-s_{13})A_1$
is the sum
\begin{align*} \sum_{\stackrel{0 \leq i < j < k}{i+j+k ~=~ 3m+1}}
\bigg( &(C_{[j,k]}~-~ C_{[i,j]})(x_1^ix_2^jx_3^k ~-~ x_1^kx_2^jx_3^i)~+~\\
&(C_{[i,k]}~-~ C_{[j,k]})(x_1^jx_2^kx_3^i ~-~ x_1^ix_2^kx_3^j) ~+~ \\
&(C_{[i,j]}~-~ C_{[i,k]})(x_1^kx_2^ix_3^j ~-~ x_1^jx_2^ix_3^k)
\bigg)
~+~\\
\sum_{\stackrel{0 \leq i, j}{2i+j ~=~ 3m+1}} \bigg( &(C_{[i,j]}
~-~ C_{[i,i]})(x_1^ix_2^ix_3^j ~-~ x_1^jx_2^ix_3^i) \bigg).
\end{align*}

\noindent We can now discover properties of the coefficients by
focusing on one summand at a time.  For instance, given a specific
composition $[i,j,k]$ of $3m+1$ such that $0 \leq i < j < k$, the
fact that $i +j + k > 3m$ means that the largest exponent, namely
$k$, will be greater than $m$.  However,  $A_1$ $m$-quasiinvariant
means that $(x_1-x_3)^{2m+1} \bigg | (1-s_{13})A_1$ and thus the
highest power of $x_2$ that can appear in $(1-s_{13})A_1$ will be
$(3m+1) - (2m+1) = m$. Thus $x_2^k$ cannot appear in a term of
$(1-s_{13})A_1$ with a nonzero coefficient, and thus we obtain
$C_{[i,k]} = C_{[j,k]}$.  If both the exponents $j$ and $k$ happen
to be greater than $m$, then by similar logic we conclude that
$C_{[i,j]} = C_{[i,k]} = C_{[j,k]}$.  Finally, if we are given the
composition $[i,i,j]$ with $i>m$, we see that $C_{[i,j]} =
C_{[i,i]}$.  We summarize these conditions here:
\begin{align}
&C_{[i,k]}=C_{[j,k]}&& \mathrm{ when~~} i<j<k &&\label{cpairs}\\
&C_{[i,j]} = C_{[i,k]} = C_{[j,k]} &&\mathrm{ when~~} i<j<k, &&j
>
m \label{ctrips}\\
&C_{[i,j]} = C_{[i,i]} &&\mathrm{ when~~} i=j, && i > m.
\label{clast}
\end{align}

\noindent The idea now will be to subtract certain symmetric
functions from $A_1$ in order to get rid of exponents of $x_2$ and
$x_3$ greater than $m$, without changing the equivalence class of
$A_1$ in the quotient. For every triplet $\{i,j,k\}$ of exponents
with $i<j<k, j \le m$, we see that $A_1 - C_{[i,k]}m_{i,j,k}$ has
\begin{align}
(C_{[i,j]} - C_{[i,k]})(x_1^kx_2^ix_3^j + x_1^kx_2^jx_3^i)
\end{align}
as the only monomials with exponent sequence a permutation of
$(i,j,k)$, by (\ref{cpairs}). (Here $m_{i,j,k}$ is the monomial
symmetric function with exponents $i,j,k$). For every triplet
$\{i,j,k\}$ of exponents with $i<j<k, j > m$, we have that $A_1 -
C_{[i,k]}m_{i,j,k}$ has no monomials with exponent sequence a
permutation of $(i,j,k)$, by (\ref{ctrips}).  For every remaining
triplet $\{i,i,j\}$ of exponents we see that $A_1 -
C_{[i,j]}m_{i,i,j}$ has
\begin{align}
(C_{[i,i]} - C_{[i,j]})x_1^jx_2^ix_3^i
\end{align}
as the only monomial with exponent sequence a permutation of
$(i,i,j)$, which by (\ref{clast}) is only nonzero when $i \leq m$.
Thus, after subtracting appropriate symmetric functions we are
left with a sum containing only monomials such that the exponents
of $x_2$ and $x_3$ are less than or equal to $m$.  This gives the
stated result for $A_1$.

\vspace{1.5em} \noindent Since $A_2$ has degree $3m+2$, the
highest power of $x_2$ that can appear in $(1-s_{13})A_2$ is
$m+1$.  Thus any composition $[i,j,k]$ such that $0 \leq i < j <
k$ and $i +j + k = 3m+2$ will have to satisfy $k > m+1$, which
will allow us to equate certain coefficients as above.  Any
composition where $0 \leq i, j$ and $2i+j = 3m+2$ will only yield
three terms, two of which have the same coefficient.  Either way,
we will analogously be able to use appropriate symmetric functions
to subtract from $A_2$ so that monomials with powers of $x_2$ or
$x_3$ exceeding $m+1$ will disappear.
\end{proof}

\noindent Later on, we will demonstrate that, for $A_2$, we can
strengthen the result of Lemma~\ref{betterform}. Namely we will
prove that there exists a quasiinvariant $A_2$ of degree $3m+2$
that satisfies the properties of Lemma~\ref{linind} and is of the
form
\begin{align} \label{3m+2deferred} A_2 = \sum_{0 \leq i \leq j \leq m}
\tilde{C}_{[i,j]} x_1^{3m+2-i-j} m_{[i,j]}(x_2,x_3).
\end{align}
Note that the indices of the sum are now less than $m+1$.  The
proof of this will require the explicit construction of $A_1$, and
will be necessary to explicitly construct $A_2$.

\section{Relations satisfied by the coefficients $C_{[i,j]}$}
\label{relations}

\noindent  In this section we show the $C_{[i,j]}$ satisfy certain
relations. We begin by setting $d=3m+1$, and

\begin{align*}  A_{i,j,k,l} =
\begin{cases}
 {i \choose k} {d-i-k \choose l} -  {i \choose k} {2i-k \choose l} &\mathrm{if~} i=j, \\
 {i \choose k} {d-j-k \choose l} + {j \choose k} {d-i-k \choose l} - \left( {i \choose k} + {j \choose k} \right) {i+j-k \choose l}
  &\mathrm{otherwise.}
\end{cases}
\end{align*}

\noindent We can now state the main result of this section.

\begin{Lem} The coefficients $C_{[i,j]}$ satisfy the linear equations
\begin{align}\label{Cij}
\sum_{0 \leq j \leq i \leq m} A_{i,j,k,l} C_{[i,j]} = 0
\end{align}
\noindent for $k \in \{0,\dots,m\}$ and $l \in
\{1,3,5,\dots,2m-1\}$.
\end{Lem}

\begin{proof} By definition, if $i > j$, then
$C_{[i,j]}$ is the coefficient of $$x_1^{d-i-j}m_{[i,j]}(x_2,x_3) =
x_1^{d-i-j} \Big ( x_2^ix_3^j + x_2^jx_3^i \Big ) .$$

\noindent If instead $i=j$, then  $C_{[i,j]}$ is the coefficient
of $x_1^{d-2i} x_2^i x_3^i$.  Consequently, inside of
$(1-s_{13})A_1$, $C_{[i,j]}$ is the coefficient of the polynomial
$$x_1^{d-i-j}x_2^ix_3^j + x_1^{d-i-j}x_2^jx_3^i -
x_1^ix_2^jx_3^{d-i-j} - x_1^jx_2^ix_3^{d-i-j}$$ if $i > j$ and
$$x_1^{d-2i} x_2^i x_3^i - x_1^i x_2^i x_3^{d-2i}$$ if $i=j$.  Using
the substitutions $~y_1 = x_2-x_1$ and $y_2 = x_1-x_3$, we rewrite
these polynomials.  For the case $i=j$ we have
\begin{align*}
(1-s_{13})A_1 \Big|_{C_{[i,i]}}
&= x_1^{d-2i}(y_1+x_1)^i x_3^i - x_1^i(y_1+x_1)^i x_3^{d-2i} \\
&= \sum_{k=0}^i {i \choose k}x_1^{d-i-k}y_1^kx_3^i - {i \choose k}
   x_1^{2i-k}y_1^{k}x_3^{d-2i} \\
&= \sum_{k=0}^i {i \choose k}(y_2+x_3)^{d-i-k}y_1^kx_3^i -
   {i \choose k}(y_2+x_3)^{2i-k}y_1^{k}x_3^{d-2i} \\
&= \sum_{k=0}^i {i \choose k}\Bigg ( \sum_{l = 0}^{d-i-k}
   {d-i-k \choose l} y_1^k y_2^{l} x_3^{d-k-l} - \sum_{l = 0}^{2i-k}
   {2i-k \choose l} y_1^{k}y_2^{l}x_3^{d-k-l} \Bigg ) \\
&= \sum_{k=0}^i \sum_{l=0}^{\max(d-i-k,~~2i-k)} {i \choose k}
   \Bigg( {d-i-k \choose l} - {2i-k \choose l} \Bigg)
   y_1^ky_2^lx_3^{d-k-l} \\
&= \sum_{k=0}^i \sum_{l=0}^{\max(d-i-k,~~2i-k)} A_{i,i,k,l}
   y_1^ky_2^lx_3^{d-k-l}
\end{align*}

\noindent For $i>j$ we have
\begin{align*}
(1-s_{13})A_1 \Big|_{C_{[i,j]}}
&= x_1^{d-i-j}(y_1+x_1)^ix_3^j + x_1^{d-i-j}(y_1+x_1)^jx_3^i \\
&- x_1^i(y_1+x_1)^jx_3^{d-i-j} - x_1^j(y_1+x_1)^ix_3^{d-i-j} \\
&= \sum_{k=0}^i {i \choose k}x_1^{d-j-k}y_1^kx_3^j +
    {j \choose k}x_1^{d-i-k}y_1^kx_3^i \\
&- {i \choose k}x_1^{i+j-k}y_1^kx_3^{d-i-j} - {j \choose k}
    x_1^{i+j-k}y_1^kx_3^{d-i-j} \\
&= \sum_{k=0}^i {i \choose k}(y_2+x_3)^{d-j-k}y_1^kx_3^j +
    {j \choose k}(y_2+x_3)^{d-i-k}y_1^kx_3^i \\
&- {i \choose k}(y_2+x_3)^{i+j-k}y_1^kx_3^{d-i-j} - {j \choose k}
    (y_2+x_3)^{i+j-k}y_1^kx_3^{d-i-j} \\
&= \sum_{k=0}^i \sum_{l=0}^{\max\{d-j-k,~i+j-k\}} \Bigg ( {i
    \choose k} {d-j-k \choose l} + {j \choose k} {d-i-k \choose l} \\
&- {i \choose k}{i+j-k \choose l} - {j \choose k}{i+j-k \choose l}
    \Bigg ) y_1^ky_2^lx_3^{d-k-l} \\
&= \sum_{k=0}^i \sum_{l=0}^{\max\{d-j-k,~i+j-k\}} A_{i,j,k,l}
    y_1^ky_2^lx_3^{d-k-l}
\end{align*}

\noindent By definition, $A_1$ is $m$-quasiinvariant if and only
if $(1-s_{13})A_1$ is divisible by $y_2^{2m+1}$. Solving the
equations implies that $(1-s_{13})A_1\bigg|_{C_{[i,j]}}$ has even
order or order greater than $2m-1$ with respect to $y_2$.  Since
$(1-s_{13})A_1$ is divisible by an odd power of $(x_1-x_3)$, we
make the following statement:  for fixed $k \in \{0,\dots,m\}$ and
fixed odd $l < 2m+1$, we must have
\begin{align*}
\sum_{0 \le j \le i \le m} A_{i,j,k,l} C_{[i,j]} y_1^k y_2^l
x_3^{d-k-l} = 0.
\end{align*}
The lemma is an immediate consequence.
\end{proof}

\section{The coefficients have a one-dimensional solution
space} \label{matrices}

Once we verify that the relations in (\ref{Cij}) have a one-dimensional solution
space, it is a straightforward (although time-intensive) process
to find a representative solution.  This will allow us to
explicitly construct $A_1$, for which we currently have only an
existence proof.  We begin by computing the determinants of
certain matrices, beginning with Theorem~\ref{Gendet}, stated in the
introduction.

\renewcommand{\proofname}{Proof of Theorem~\ref{Gendet}}
\begin{proof}

\noindent We first show that the number of lattice paths from
$(c+jd, c+jd)$ to $(0,a+ib)$ which avoid the line $y=-x +(c+e)$ is
$\binom{a+bi}{c+dj} - \binom{a+bi}{e-dj}$. Consider the following
two diagrams:

\vspace{1.5em}
\begin{figure}[hbpt]
\centerline{\input{jan2}} \caption{Counting paths from $(c+jd, c+jd)$ to
$(0,a+ib)$.\label{f-jan2}}
\end{figure}
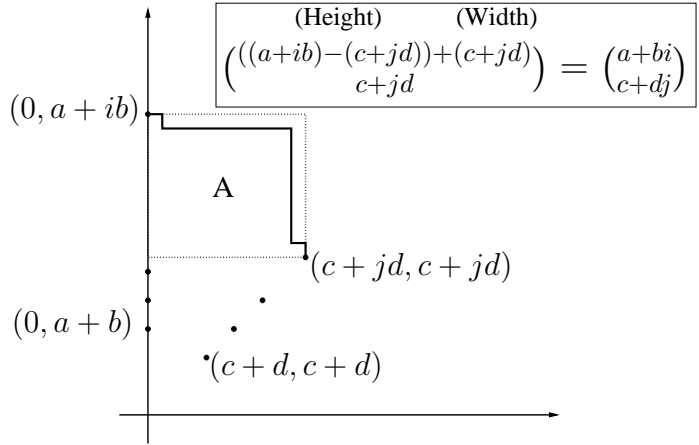

\begin{figure}[hbpt]
\centerline{\input{jan3}} \caption{`Bad' paths from $(c+jd, c+jd)$ to $(0,a+ib)$.\label{f-jan3}}
\end{figure}
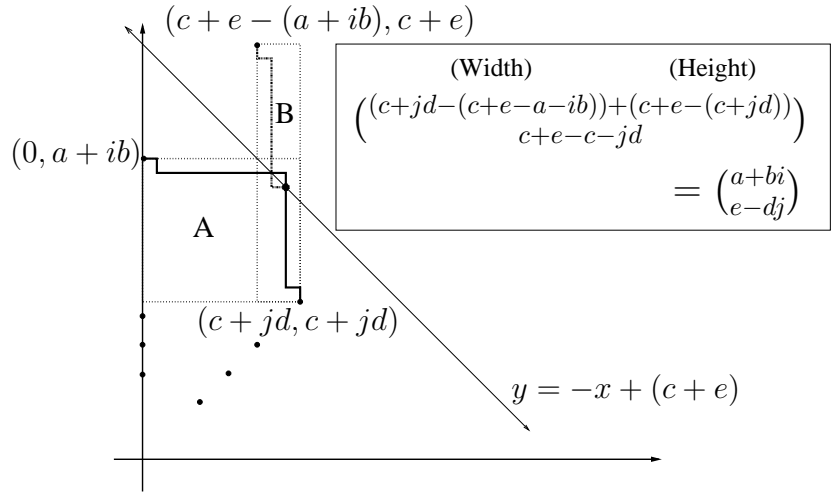

\vspace{1.5em}\noindent The number of bad paths in rectangle $A$,
namely the ones that go through the forbidden line, is in
bijection with the number of total paths in rectangle $B$; we
replace WEST steps with NORTH steps and NORTH steps with WEST
steps following the first touch of the forbidden line.  This is
known as Andr\'{e}'s Reflection Principle \cite{Andre}. Thus the
number of good paths in rectangle $A$ is exactly the correct
difference of binomials.

\vspace{1.5em}\noindent This shown, a classical involution of
Lindstr\"{o}m \cite{Lind} and Gessel-Viennot \cite{Gessel} shows
that when the entries of a matrix count paths, the determinant
counts families of non-intersecting paths.  This completes the
proof.
\end{proof}
\renewcommand{\proofname}{Proof}

\noindent We now are in a position to prove
Theorem~\ref{specthm}, as stated in the introduction.
\renewcommand{\proofname}{Proof of Theorem~\ref{specthm}}
\begin{proof}

\noindent We begin by considering a more general form of this
matrix and factoring it. This factorization was suggested by an
argument of Gessel and Viennot \cite{Gessel}:
\[
\det \Bigg| \binom{a_i}{b_j} - \binom{c-a_i}{b_j} \Bigg|_{i,j=1}^k
= \]
\[ \det \Bigg| \frac{\binom{c}{b_{k-i+1}}}{\binom{c}{a_{k-j+1}}} \cdot
\left( \binom{c-b_{k-i+1}}{c-a_{k-j+1}} -
\binom{c-b_{k-i+1}}{a_{k-j+1}} \right) \Bigg|_{i,j=1}^k =\]
\[ \frac{\binom{c}{b_1} \dots \binom{c}{b_k}} {\binom{c}{a_1} \dots
\binom{c}{a_k}} \cdot \det \Bigg| \binom{c-b_{k-i+1}}{c-a_{k-j+1}}
- \binom{c-b_{k-i+1}}{a_{k-j+1}} \Bigg|_{i,j=1}^k
\]

\noindent Proposition 14 of \cite{Gessel} used an analogous
factorization for the determinant of a matrix of single binomial
coefficients. Our factorization also works by the symmetry ${c
\choose a_i} = {c \choose c-a_i}$.  This implies that the same
quotient of binomials can be factored out of both terms that
appear as a difference in our entries.  Returning to the proof of
Theorem \ref{specthm}, we let $a_i = C + \alpha i$, $b_j = E +
\beta j$, and $c = C+D$ and find

\begin{align*}
\det \Bigg| \binom{C+\alpha i}{E+\beta j} &- \binom{D-\alpha
i}{E+\beta j} \Bigg|_{i,j=1}^k =
\frac{\binom{C+D}{E+\beta}\binom{C+D}{E+2\beta} \dots
\binom{C+D}{E+n\beta}}
{\binom{C+D}{C+\alpha}\binom{C+D}{C+2\alpha}
\dots \binom{C+D}{C+n\alpha}} ~~\bullet \\
&\det \Bigg| \binom{C+D-E-(k-i+1)\beta}{D - (k-j+1)\alpha} -
\binom{C+D-E-(k-i+1)\beta}{C + (k-j+1)\alpha} \Bigg|_{i,j=1}^k.
\end{align*}

\vspace{1.5em}\noindent Notice that now the tops of the binomial
coefficients are the same and the bottoms are different. This
allows us to apply Theorem \ref{Gendet} to obtain the result.
\end{proof}
\renewcommand{\proofname}{Proof}

\noindent We now see how these results can help us with our system
of equations.  Notice that in (\ref{Cij}) there are ${m + 2
\choose 2}$ coefficients $C_{[i,j]}$ and $m(m+1)$ equations.  We
define $B_m$ as the restriction of the matrix given by (\ref{Cij})
to the $\big( {m+2 \choose 2} - 1\big) \times \big( {m+2 \choose
2} - 1\big)$ sub-matrix where $[i,j] \not = [m,m]$, $0 \leq k \leq
m-1$ and $l \in \{2m-2k-1,\dots,2m-3,2m-1\}$ or $k=m$ and $l \in
\{1,3,5,\dots, 2m-1\}$.

\begin{Lem}The matrix $B_m$ is nonsingular
\end{Lem}
\begin{proof}
By using an ordering for the pairs $(k,l)$ where the $k$'s
increase and the $l$'s decrease while lexicographically ordering
the $[i,j]$s, the matrix $B_m$ becomes block triangular.
Furthermore, there is one block of size $1$, one block of size
$2$, $\dots$, one block of size $m-1$, and two blocks of size $m$.
This block triangularity follows from the fact that for $i,j$ such
that $0 \leq j \leq i < k$, then ${i \choose k} = {j \choose k} =
0$ and thus the $A_{i,j,k,l}$'s of equation (\ref{Cij}) are all
zero.

\vspace{1.5em} \noindent Furthermore, the entries of $B_m$
inside these blocks, where $j$ runs over the interval $0 \leq
j \leq i = k$, are much simpler than the general case.  For such
$i,j$'s, the $A_{i,j,k,l}$'s of equation (\ref{Cij}) simplify to

\begin{align*}  A_{k,j,k,l} =
\begin{cases}
 {k \choose k} {d-2k \choose l} -  {k \choose k} {k \choose l} &\mathrm{if~} j= k, \\
 {k \choose k} {d-j-k \choose l} + {j \choose k} {d-2k \choose l} - \left( {k \choose k} + {j \choose k} \right) {j \choose l}
  &\mathrm{otherwise.}
\end{cases}
\end{align*}

\noindent But since ${k \choose k} = 1$ and ${j \choose k} = 0$ if
$j < k$ we obtain

\begin{align} \label{Bsimplified} A_{k,j,k,l} =
 {d- j - k \choose l} -  {k \choose l}.
\end{align}

\noindent For $f \in \{1,2,\dots, m\}$, we let $B^{f,m}$ denote
the $f^{th}$ block matrix on the diagonal of $B_m$, which forces
$f=k+1$, and set $B^m$ to be the final block matrix. Setting
$d=3m+1$ and utilizing (\ref{Bsimplified}) allows us to describe
the entries of these blocks as follows:

\vspace{1.5em}\noindent For $j \in \{ 0,\dots, f-1 \}$ and $l \in
\{2m-2f+1,2m-2f-1,\dots, 2m-1\}$,
\begin{align*}
B^{f,m}_{l,j} = {3m+1-j-(f-1) \choose l} - {j \choose l}
\end{align*}
and for $j \in \{0,\dots, m-1\}$ and $l \in \{1,3,\dots, 2m-1\}$,
\begin{align*}
B^{m}_{l,j} = {2m+1-j \choose l} - {j \choose l}.
\end{align*}

\noindent At this point, we re-index the matrix $B^{f,m}$,
replacing the current indices of $j$ and $l$ with the standard
indices $i,j \in \{1,\dots,f\}$.  This gives

\begin{align} \label{block1}
B^{f,m} = \Bigg| \binom{3m+1-(j-1)-(f-1)}{2m+1-2i} -
\binom{j-1}{2m+1-2i} \Bigg|_{i,j=1}^f
\end{align}

\noindent and

\begin{align} \label{block2}
B^{m} = \Bigg| \binom{2m+1-(j-1)}{2m+1-2i} - \binom{j-1}{2m+1-2i}
\Bigg|_{i,j=1}^m.
\end{align}

\noindent Applying Theorem~\ref{specthm} to the transpose of this
matrix, we find the determinant of (\ref{block1}) is

\[ \frac{\binom{3m+2-f}{2m-1} \binom{3m+2-f}{2m-3} \cdots
\binom{3m+2-f}{2m-2f+1}} {\binom{3m+2-f}{3m+2-f}
\binom{3m+2-f}{3m+1-f} \cdots \binom{3m+2-f}{3m-2f+3}} \cdot |
\mathcal{F} |
\]
where $\mathcal{F}$ is the set of families of non-intersecting
lattice paths from $\{(0,0),(1,1),\dots,(f-1,f-1)\}$ to
$\{(0,m-f+3),(0,m-f+5),\dots,(0,m+f+1)\}$ which stay below the
line $y=-x+3m+2-f$.  Since this family of paths is non-empty, we
conclude that the matrices $B^{f,m}$ are non-singular for $f \in
\{1,2,\dots, m\}$.  Similarly we find that the determinant of
(\ref{block2}) is positive and thus $B^{m}$ is also non-singular.
Since the diagonal blocks of $B_m$ are non-singular, the matrix
$B_m$ must also be.
\end{proof}

\noindent An example may help to clarify things at this point.
When $m=3$ and $d=10$, we have the matrix
$$\begin{bmatrix}
252& 378 & 126 & 308& 182 & 56 & 273 & 147 & 75  \\
0  & 126 & 56  & 252& 133 & 42 & 378 & 174 & 75  \\
0  & 84  & 56  & 168& 147 & 68 & 252 & 184 & 125 \\
0  & 0   & 0   & 56 & 21  & 6  & 168 & 63  & 19  \\
0  & 0   & 0   & 56 & 35  & 20 & 168 & 105 & 66  \\
0  & 0   & 0   & 8  & 6   & 4  & 21  & 15  & 11  \\
0  & 0   & 0   & 0  & 0   & 0  & 21  & 6   & 1   \\
0  & 0   & 0   & 0  & 0   & 0  & 35  & 20  & 10  \\
0  & 0   & 0   & 0  & 0   & 0  & 7   & 5   & 3
\end{bmatrix}$$

\vspace{1em} \noindent This matrix is the matrix of coefficients
$A_{i,j,k,l}$ where the columns are indexed by the $[i,j]$'s and
the rows are indexed by the pairs $(k,l)$. In this example, the
columns have the order

$$[i,j] =  [0,0], ~[1,0], ~[1,1], ~[2,0], ~[2,1], ~[2,2], ~[3,0], ~[3,1],
~[3,2] $$

\noindent and the rows have the order:

$$(k,l) = (0,5), ~(1,5), ~(1,3), ~(2,5), ~(2,3), ~(2,1), ~(3,5), ~(3,3),
~(3,1). $$

\noindent We also have the following block sub-matrices:

$$B^{1,3} =\begin{bmatrix}
252
\end{bmatrix}$$

$$B^{2,3} = \begin{bmatrix}
126 & 56    \\
84  & 56
\end{bmatrix}$$

$$B^{3,3} = \begin{bmatrix}
56 & 21 &6      \\
56 & 35 &20    \\
8  & 6  &4
\end{bmatrix}$$

$$B^{3} = \begin{bmatrix}
21  & 6   & 1   \\
35  & 20  & 10  \\
7   & 5   & 3
\end{bmatrix}$$

\begin{Lem} The equations given in (\ref{Cij}) have a solution that is
unique up to scalar multiples.
\end{Lem}

\begin{proof}
Since the system in (\ref{Cij}) has an $\big( {m+2 \choose 2} -
1\big) \times \big( {m+2 \choose 2} - 1\big)$ nonsingular
sub-matrix, it must be true that the rank of the system in
(\ref{Cij}) is $\geq { m + 2 \choose 2}-1$.  Thus the null space
has dimension $\leq 1$. However, since we know by Lemma
\ref{betterform} that $A_1$ is a solution, the dimension of the
null space must be exactly one.
\end{proof}

\section{Constructing $A_2$ and a basis for the quotient}

\noindent In the first section, we showed the existence of nonzero
(in the quotient) $m$-quasiinvariants $A_1$, $A_2$ of degrees
$3m+1$ and $3m+2$, respectively, that are both symmetric with
respect to $s_{23}$. In the previous two sections we illustrated
an explicit construction of the element $A_1$.  We now give an
explicit construction of the element $A_2$, which will be linearly
independent of $A_1$. We have deferred this construction until now
since this argument is dependent on the explicit form of $A_1$. We
begin by strengthening Lemma \ref{betterform}.

\begin{Lem} \label{stronger}
There exists an $m$-quasiinvariant of degree $3m+2$, satisfying
the conditions of Lemma \ref{linind}, which has the form given in
equation (\ref{3m+2deferred}).
\end{Lem}

\begin{proof}
\noindent First, we observe that the Hilbert series (\ref{hilb})
and Lemma \ref{betterform} tell us there is a nonzero
$m$-quasiinvariant of degree $3m+1$, $$A_1 = \sum_{0 \leq i \leq j
\leq m} C_{[i,j]} x_1^{3m+1-i-j} m_{[i,j]}(x_2,x_3),$$ as well as
a nonzero $m$-quasiinvariant of degree $3m+2$, $$A_2 = \sum_{0
\leq i \leq j \leq m+1} \tilde{C}_{[i,j]} x_1^{3m+2-i-j}
m_{[i,j]}(x_2,x_3).$$

\noindent We proved in the last section that the set of possible
coefficient vectors $\langle C_{[i,j]}\rangle$ comprises a
$1$-dimensional space. Eliminating the last column of the matrix
of entries $A_{i,j,k,l}$'s is like setting the coefficient
$C_{[m,m]} = 0$.  Since the sub-matrix $B_m$ also lacks that
column and is nonsingular we conclude that the nonzero
$m$-quasiinvariant $A_1$ satisfies $C_{[m,m]} \not = 0$.
Consequently,  $e_1A_1$ has a nonzero multiple of
$x_1^{m+1}x_2^{m+1}x_3^m$ as one of its terms while at the same
time the term $x_1^mx_2^{m+1}x_3^{m+1}$ will not appear.  With no
cancellation therefore possible, the quantity
 $(1-s_{13})e_1A_1$ will contain the term $C(x_1-x_3)^{2m+1}x_2^{m+1}$
 for some nonzero $C$.

\vspace{1.5em}\noindent Since $(1-s_{13})A_2 =
(x_1-x_3)^{2m+1}(C^{\prime}x_2^{m+1} +
\mathrm{~terms~of~lower~order})$, we find that
$(1-s_{13})(C^{\prime}e_1A_1 - CA_2)$ contains no term with
$x_2^{m+1}$.  We thus re-define $A_2$ as the quantity
$C^{\prime}e_1A_1 - CA_2$ (which still meets the conditions of
Lemma~\ref{linind}). Recall that in the proof of Lemma
\ref{betterform}, the crucial step that proved the result for
$A_1$ was the fact that we could eliminate every term in
$(1-s_{13})A_1$ containing a power of $x_2$ exceeding $m$.  Now we
can utilize this fact for $A_2$ also. The rest of the proof goes
through as before and we conclude that $A_2$ can be written as

\begin{align} \label{Ctwiddle}
\sum_{0 \leq i \leq j \leq m} \tilde{C}_{[i,j]} x_1^{3m+2-i-j}
m_{[i,j]}(x_2,x_3).
\end{align}
\end{proof}

\noindent We now examine how the construction of $A_1$ can be
applied to construct $A_2$. In section \ref{relations}, we used
the fact that $A_1$ had the form
$$\sum_{0 \leq i \leq j \leq m} C_{[i,j]} x_1^{3m+1-i-j}
m_{[i,j]}(x_2,x_3).$$ \noindent to obtain a linear system of
relations that the $C_{[i,j]}$'s satisfy.  Since we now know that
$A_2$ has an analogous form, namely (\ref{Ctwiddle}), we can apply
the same proof (setting $d = 3m+2$) to obtain an analogous system for
the $\tilde{C}_{[i,j]}$'s.

\vspace{1.5em}\noindent These coefficients can be explicitly
computed by finding the null space of the matrix given by the
linear system

\begin{align} \label{Ctildeij}
\sum_{0 \leq j \leq i \leq m} A_{i,j,k,l} \tilde{C}_{[i,j]} = 0
\end{align}

\noindent for $k \in \{0,\dots,m\}$ and $l \in
\{1,3,5,\dots,2m-1\}$.  As in the $A_1$ case, this null space is
$1$-dimensional and we prove this by showing that the matrix
$\tilde{B}_m$ is nonsingular, where $\tilde{B}_m$ is the
restriction of the matrix given by (\ref{Ctildeij}) to the $\big(
{m+2 \choose 2} - 1\big) \times \big( {m+2 \choose 2} - 1\big)$
sub-matrix where $[i,j] \not = [m,m]$, $0 \leq k \leq m-1$ and $l
\in \{2m-2k-1,\dots,2m-3,2m-1\}$ or $k=m$ and $l \in
\{1,3,5,\dots, 2m-1\}$.

\vspace{1.5em}\noindent The matrix $\tilde{B}_m$ is block
triangular and thus we prove that it is nonsingular by proving
that its blocks

\begin{align} \label{block3}
\tilde{B}^{f,m} = \Bigg| \binom{3m+2-(j-1)-(f-1)}{2m+1-2i} -
\binom{j-1}{2m+1-2i} \Bigg|_{i,j=1}^f
\end{align}

\noindent for $f \in \{1, 2, \dots, m\}$ as well the additional
block
\begin{align} \label{block4}
\tilde{B}^{m} = \Bigg| \binom{2m+2-(j-1)}{2m+1-2i} -
\binom{j-1}{2m+1-2i} \Bigg|_{i,j=1}^m
\end{align}

\noindent are nonsingular.  We proceed identically to our
computation of the determinant of (\ref{block1}). We find that the
determinant of (\ref{block3}) is a positive scalar multiplied by
the number of families of non-intersecting lattice paths from
$\{(0,0),(1,1),\dots,(f-1,f-1)\}$ to
$\{(0,m-f+4),(0,m-f+6),\dots,(0,m+f+2)\}$ which stay below the
line $y=-x+3m+3-f$.  Since such paths exist, this determinant is
positive.  Similarly we find that the determinant of
(\ref{block4}) is positive and thus our construction of $A_2$ is
valid.

\begin{Thm} The set $\{1, A_1, s_{12}(A_1), A_2, s_{12}(A_2),
\Delta^{2m+1}(x)\}$ is a basis for the quotient $QI_m \big/
\big\langle (e_1,e_2,e_3) \big\rangle$.
\end{Thm}

\begin{proof} It remains only to prove the independence of $\{A_1,
s_{12}(A_1), A_2, s_{12}(A_2) \}$ in the quotient.  By examining
the Hilbert series of $QI_m$ (\ref{hilb2}), we find that the
subspace of $\q[x_1,x_2,x_3]$ consisting of $3m+2$ dimensional
$m$-quasiinvariants which are not symmetric is 4 dimensional. Thus
it is spanned by $e_1A_1, e_1(s_{12})A_1$ and two other elements.
Since we have shown that $A_i$ and $s_{12}A_i$ are linearly
independent for $i \in \{1,2\}$, it remains to show that there is
no nontrivial collection of constants $c_1, c_2, c_3, c_4$ such
that
\begin{align}
\label{lindep}  c_1e_1A_1 + c_2e_1s_{12}A_1 + c_3A_2 +
c_4s_{12}A_2 &= 0.
\end{align}
We first note that
\begin{align}\label{e1a1}
A_2 \ne c e_1 A_1\mathrm{~~for~any~~}c.
\end{align}
This is seen by examining the terms containing $x_2^{m+1}$ in
each, as was done in the proof of Lemma \ref{stronger}. Now assume
that (\ref{lindep}) held.  Applying $s_{13}\pi_2$ gives
\begin{align}
c_2e_1A_1 + c_4A_2 &= 0
\end{align}
which is in immediate contradiction of (\ref{e1a1}), unless
$c_2=c_4=0$.  Returning to (\ref{lindep}) gives
\begin{align}
c_1e_1A_1 + c_3A_2 &= 0.
\end{align}
Again (\ref{e1a1}) forces $c_1=c_3=0$. This completes the proof.
\end{proof}

\noindent We have thus reduced the problem of finding a basis for
the quasiinvariants of $S_3$ to finding the 1-dimensional
nullspace of particular matrices, a computation easily carried out
by computer. We have used this technique to explicitly compute the
basis for several small values of $m$. We conclude with the
following examples:

\vspace{1.5em} \noindent For $m = 1$,
\begin{align*}
A_{1} &= x_1^4 - 2x_1^3(x_2+x_3) + 6x_1^2(x_2x_3) \\
A_{2} &= x_1^5 - {5 \over 3}x_1^4(x_2+x_3) + {10 \over
3}x_1^3(x_2x_3).
\end{align*}

\noindent For $m = 2$, {
\begin{align*}
A_1 &= x_1^7 - {7 \over 2}x_1^6(x_2+x_3) + 14x_1^5(x_2x_3) + {7
    \over 2} x_1^5(x_2^2+x_3^2) \\
&- {35 \over 2}x_1^4(x_2^2x_3 + x_2x_3^2) + 35x_1^3x_2^2x_3^2  \\
A_2 &= x_1^8 - {16 \over 5}x_1^7(x_2+x_3) + {56 \over
    5}x_1^6(x_2x_3) + {14 \over 5} x_1^6(x_2^2+x_3^2)\\
&- {56 \over 5}x_1^5(x_2^2x_3 + x_2x_3^2) + 14x_1^4x_2^2x_3^2 .
\end{align*} }

\newpage
\noindent \bf Acknowledgements. \rm We would like to thank
Christian Krattenthaler for his helpful suggestions and his
wonderful reference for finding determinants \cite{Kratt}. We are
also indebted to Adriano Garsia for introducing us to this
subject. We are thankful for his guidance and support during this
project.

\end{document}

%% file: jan2.tex
\begin{picture}(0,0)%
\includegraphics{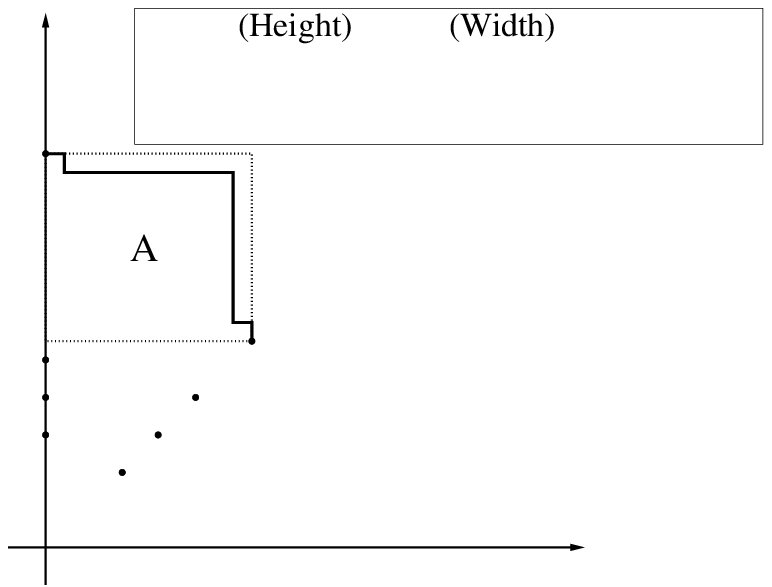}%
\end{picture}%
\setlength{\unitlength}{1184sp}%
\begingroup\makeatletter\ifx\SetFigFont\undefined%
\gdef\SetFigFont#1#2#3#4#5{%
  \reset@font\fontsize{#1}{#2pt}%
  \fontfamily{#3}\fontseries{#4}\fontshape{#5}%
  \selectfont}%
\fi\endgroup%
\begin{picture}(14412,9270)(-1424,-9394)
\put(3001,-1711){\makebox(0,0)[lb]{\smash{\SetFigFont{5}{6.0}{\rmdefault}{\mddefault}{\updefault}{\color[rgb]{0,0,0}{\Large $\binom{\left((a+ib)-(c+jd)\right)+(c+jd)}{c+jd} =\binom{a+bi}{c+dj} $}}%
}}}
\put(-1349,-7036){\makebox(0,0)[lb]{\smash{\SetFigFont{5}{6.0}{\rmdefault}{\mddefault}{\updefault}{\color[rgb]{0,0,0}{\large $(0,a+b)$}}%
}}}
\put(-1424,-2536){\makebox(0,0)[lb]{\smash{\SetFigFont{5}{6.0}{\rmdefault}{\mddefault}{\updefault}{\color[rgb]{0,0,0}{\large $(0,a+ib)$}}%
}}}
\put(4876,-5836){\makebox(0,0)[lb]{\smash{\SetFigFont{5}{6.0}{\rmdefault}{\mddefault}{\updefault}{\color[rgb]{0,0,0}{\large $(c+jd,c+jd)$}}%
}}}
\put(2776,-7936){\makebox(0,0)[lb]{\smash{\SetFigFont{5}{6.0}{\rmdefault}{\mddefault}{\updefault}{\color[rgb]{0,0,0}{\large $(c+d,c+d)$}}%
}}}
\end{picture}

%% file: jan3.tex
\begin{picture}(0,0)%
\includegraphics{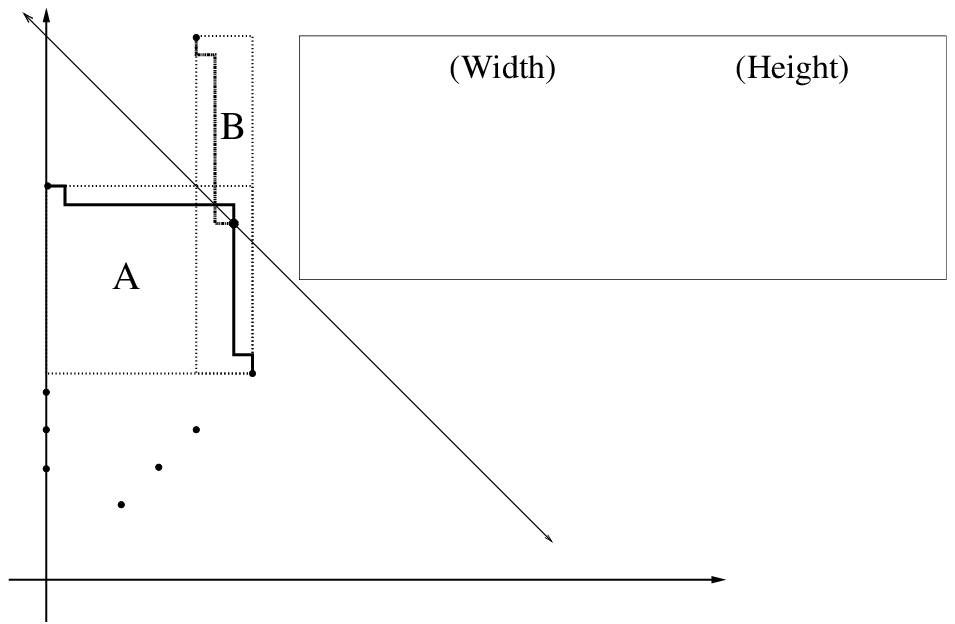}%
\end{picture}%
\setlength{\unitlength}{1184sp}%
\begingroup\makeatletter\ifx\SetFigFont\undefined%
\gdef\SetFigFont#1#2#3#4#5{%
  \reset@font\fontsize{#1}{#2pt}%
  \fontfamily{#3}\fontseries{#4}\fontshape{#5}%
  \selectfont}%
\fi\endgroup%
\begin{picture}(17187,9891)(-4274,-9394)
\put(-1049,314){\makebox(0,0)[lb]{\smash{\SetFigFont{5}{6.0}{\rmdefault}{\mddefault}{\updefault}{\color[rgb]{0,0,0}{\large $(c+e-(a+ib),c+e)$}}%
}}}
\put(-374,-5986){\makebox(0,0)[lb]{\smash{\SetFigFont{5}{6.0}{\rmdefault}{\mddefault}{\updefault}{\color[rgb]{0,0,0}{\large $(c+jd,c+jd)$}}%
}}}
\put(9601,-3286){\makebox(0,0)[lb]{\smash{\SetFigFont{5}{6.0}{\rmdefault}{\mddefault}{\updefault}{\color[rgb]{0,0,0}{\Large $= \binom{a+bi}{e-dj}$}}%
}}}
\put(2851,-1786){\makebox(0,0)[lb]{\smash{\SetFigFont{5}{6.0}{\rmdefault}{\mddefault}{\updefault}{\color[rgb]{0,0,0}{\Large $\binom{\left(c+jd-(c+e-a-ib)\right)+\left(c+e-(c+jd)\right)}{c+e-c-jd}$}}%
}}}
\put(6226,-7411){\makebox(0,0)[lb]{\smash{\SetFigFont{5}{6.0}{\rmdefault}{\mddefault}{\updefault}{\color[rgb]{0,0,0}{\large $y=-x+(c+e)$}}%
}}}
\put(-4274,-2461){\makebox(0,0)[lb]{\smash{\SetFigFont{5}{6.0}{\rmdefault}{\mddefault}{\updefault}{\color[rgb]{0,0,0}{\large $(0,a+ib)$}}%
}}}
\end{picture}

%% file: 9-13-040-Quasi.bbl
\vspace{2em}


\begin{thebibliography}{99}

\bibitem{FV} M. Feigin and A. P. Veselov, Quasiinvariants of Coxeter groups and m-harmonic
polynomials, Intern. Math. Res. Notices, 2002, No. 10, 521-545.

\bibitem{FV2} G. Felder and A. P. Veselov, Action of Coxeter Groups on m-Harmonic polynomials and KZ equations, math.QA/0108012

\bibitem{EG} P. Etingof and V. Ginzburg, On m-quasiinvariants of a Coxeter group, math.QA/0106175

\bibitem{Andre} W. Feller, An introdction to probability theory
and its applications, Vol. I (John Wiley \& Sons, 1968).

\bibitem{Gessel} I. Gessel and G. Viennot, Binomial determinants, paths and hook length formulae, Adv. in Math., 58 (1985), 300-321.

\bibitem{Kane} R. Kane, Reflection Groups and Invariant Theory (Springer-Verlag New York, inc. 2001).

\bibitem{Kratt} C. Krattenthaler, Advanced determinant calculus, S\'{e}minaire Lotharingien Combin. 42
(``The Andrews Festschrift") (1999), Article B42q, 67 pp.

\bibitem{Lind} B. Lindstr\"{o}m, On the verctor representations
of induced matroids, Bull. London Math Soc., 5 (1973), 85-90.

\bibitem{EC2} R.P. Stanley, Enumerative Combinatorics, Vol. II (Cambridge University Press, 1999).

\end{thebibliography}
